\documentclass{amsart}

\newtheorem{theorem}{Theorem}
\newcommand{\bt}{\begin{theorem}}
\newcommand{\et}{\end{theorem}}

\newtheorem{lemma}{Lemma}
\newcommand{\bl}{\begin{lemma}}
\newcommand{\el}{\end{lemma}}

\newtheorem{corollary}{Corollary}
\newcommand{\bc}{\begin{corollary}}
\newcommand{\ec}{\end{corollary}}

\newtheorem{problem}{Problem}
\newcommand{\bp}{\begin{problem}}
\newcommand{\ep}{\end{problem}}

\newcommand{\beq}{\begin{equation}}
\newcommand{\eeq}{\end{equation}}
\newcommand{\benum}{\begin{enumerate}}
\newcommand{\eenum}{\end{enumerate}}

\newcommand{\Z}{\ensuremath{\mathbf Z}}

\newcommand{\N}{\ensuremath{ \mathbf N }}

\newcommand{\R}{\ensuremath{ \mathbf{R} }}
\newcommand{\card}{\ensuremath{\text{card}}}

\newcommand{\mck}{\ensuremath{\mathcal K}}
\newcommand{\mcs}{\ensuremath{\mathcal S}}
\newcommand{\mct}{\ensuremath{\mathcal T}}
\newcommand{\mcu}{\ensuremath{\mathcal U}}
\newcommand{\mcx}{\ensuremath{\mathcal X}}

\begin{document}

\title{Families of linear semigroups with intermediate growth}
\subjclass[2000]{Primary 20M05,20F69,20F65,11B75.}
\keywords{Semigroups, intermediate growth, semigroup growth, group growth, sumsets, Sidon sets, $\hat{B}_h$-sequences, partition asymptotics, additive number theory.}

\author{Melvyn B. Nathanson}
\address{Lehman College (CUNY),Bronx, New York 10468}
\email{melvyn.nathanson@lehman.cuny.edu}

\thanks{This work was supported in part by grants from the NSA Mathematical Sciences Program and the PSC-CUNY Research Award Program.}

\date{\today}

\maketitle

\begin{abstract}
Methods from additive number theory are applied to construct families of finitely generated linear semigroups with intermediate growth.
\end{abstract}

\section{Growth functions of finitely generated semigroups}
Let $\mathcal{S}$ be a finitely generated semigroup and let $A$ be a set of generators for $\mathcal{S}$.  Every element $x \in \mathcal{S}$ can be written as a word with letters from the set $A$, that is, as a finite product of elements of $A$.  The \emph{length} of $x$ with respect to $A$, denoted $\ell_A(x)$, is the 
number of letters in the shortest word that represents $x$.  Note that $\ell_A(x) = 1$ if and only if $x \in A$.  We shall assume that  \mcs\ contains an identity element 1 and that $1 \notin A$.  We define $\ell_A(1) = 0.$

Let $\N = \{1,2,3,\ldots\}$ denote the set of positive integers and $\N_0 = \N \cup \{0\}$ the set of nonnegative integers.
For every  $n \in \N_0$, let $\lambda_A(n)$ denote the number of elements of \mcs\ of length exactly $n$.
We define the \emph{growth function} $\gamma_{A}(n)$ of $\mathcal{S}$ with respect to $A$ by
\[
\gamma_{A}(n) = \card(\{x\in \mathcal{S} : \ell_A(x) \leq n\}) = \sum_{m=0}^n \lambda_A(m).
\]
The function $\gamma_A(n)$ is an increasing function that counts the number of elements of \mcs\ of length at most $n$.  If $\card(A) = k $,  then, for all nonnegative integers $n$,
\[
\lambda_A(n) \leq k^n
\]
and
\[
\gamma_A(n) \leq 
\begin{cases}
n+1 & \text{if $k=1$}\\
\frac{k^{n+1}-1}{k-1}  & \text{if $k>1$.}
\end{cases}
\]

The semigroup $\mathcal{S}$ has \emph{polynomial growth} with respect to the generating set $A$ if there exist positive numbers $c$ and $d$ such that $\gamma_A(n) \leq cn^d$ for all sufficently large $n$.  In this case,  $\log \gamma_A(n) \leq \log c + d\log n$ and so
\[
\limsup_{n\rightarrow\infty} \frac{\log \gamma_A(n)}{\log n} < \infty.
\]
If $\mathcal{S}$ does not have polynomial growth, then 
\[
\limsup_{n\rightarrow\infty} \frac{\log \gamma_A(n)}{\log n} = \infty.
\]
We say that $\mathcal{S}$ has \emph{superpolynomial growth} if 
\[
\lim_{n\rightarrow\infty} \frac{\log \gamma_A(n)}{\log n} = \infty.
\]

The semigroup $\mathcal{S}$ has \emph{exponential growth} with respect to the generating set $A$ if there exists $\theta > 1$ such that $\gamma_A(n) \geq \theta^n$ for all sufficently large $n$.  In this case,  $\log \gamma_A(n) \geq n \log \theta$ and so
\[
\liminf_{n\rightarrow\infty} \frac{\log \gamma_A(n)}{n} > 0.
\]
If \mcs\ does not have exponential growth, then 
\[
\liminf_{n\rightarrow\infty} \frac{\log \gamma_A(n)}{n} = 0.
\]
We say that the $\mathcal{S}$ has \emph{subexponential growth} if 
\[
\lim_{n\rightarrow\infty} \frac{\log \gamma_A(n)}{n} = 0.
\]
The semigroup $\mathcal{S}$ has \emph{intermediate growth} with respect to $A$ if the growth of $\mathcal{S}$ is both superpolynomial and subexponential.

We review some standard facts about growth functions.
Let $A$ and $B$ be finite generating sets for $\mathcal{S}$, and let
\[
c_A = \max(\ell_A(b) : b \in B )
\]
and
\[
c_B = \max(\ell_B(a) : a \in A ).
\]
Then
\[
\gamma_{B}(n) \leq \gamma_{A}(c_An)
\]
and
\[
\gamma_{A}(n) \leq \gamma_{B}(c_B n).
\]
These inequalities imply that the growth function $\gamma_A(n)$ is polynomial, superpolynomial, subexponential, or exponential if and only if 
the growth function $\gamma_B(n)$ is, respectively,  polynomial, superpolynomial, subexponential, or exponential.

\bl    \label{intermed:lemma:fone-to-one}
Let $\mathcal{S}$ and $\mathcal{T}$ be finitely generated semigroups and let $f: S \rightarrow T$ be an injective semigroup homomorphism.  If $A$ is any finite generating set for \mcs, then there is a finite generating set $B$ for \mct\  such that 
\[
\gamma^{(S)}_A(n) \leq \gamma^{(T)}_B(n)
\]
for all $n \in \N_0.$  If the growth of $\mathcal{T}$ is polynomial or subexponential, then the growth of $\mathcal{S}$ is, respectively, polynomial or subexponential.
If the growth of $\mathcal{S}$ is superpolynomial or exponential, then the growth of $\mathcal{T}$ is, respectively, superpolynomial or exponential.
\el

\begin{proof}
Let $A$ be a finite generating set for $\mathcal{S}$, let $B'$ be any finite generating set for $\mathcal{T}$, and let $B = B' \cup f(A).$    If $\ell_A(s) = m$, then $s=a_1\cdots a_m$ for some $a_1,\ldots,a_m \in A$.  Since $f(s) = f(a_1)\cdots f(a_m)$, it follows that  $\ell_B(f(s)) \leq m.$  Since $f$ is one-to-one, it follows that $\gamma^{(S)}_A(n) \leq \gamma^{(T)}_B(n)$.  This inequality implies the statements about growth rates.
\end{proof}

\bl    \label{intermed:lemma:subsemi}
Let $\mathcal{S}$ be a finitely generated subsemigroup of a finitely generated semigroup $\mathcal{T}$.  If the growth of $\mathcal{T}$ is polynomial or subexponential, then the growth of $\mathcal{S}$ is, respectively, polynomial or subexponential.
If the growth of $\mathcal{S}$ is superpolynomial or exponential, then the growth of $\mathcal{T}$ is, respectively, superpolynomial or exponential.
\el

\begin{proof}
This follows immediately from Lemma~\ref{intermed:lemma:fone-to-one}.
\end{proof}

\bl    \label{intermed:lemma:fonto}
Let $\mathcal{S}$ and $\mathcal{T}$ be finitely generated semigroups and let $f: S \rightarrow T$ be a surjective semigroup homomorphism.  If $A$ is any finite generating set for \mcs, then $B = \{f(a): a\in A\}$  is a finite generating set for \mct, and 
\[
\gamma^{(S)}_A(n) \geq \gamma^{(T)}_B(n) 
\]
for all $n \in \N_0$.
If the growth of $\mathcal{S}$ is polynomial or subexponential, then the growth of $\mathcal{T}$ is, respectively, polynomial or subexponential.
If the growth of $\mathcal{T}$ is superpolynomial or exponential, then the growth of $\mathcal{S}$ is, respectively, superpolynomial or exponential.
\el

\begin{proof}
Let $y \in \mct.$  Since $f$ is onto, there exists $x \in \mcs$ such that $f(x) = y.$  If $\ell_A^{(S)}(x) = m,$ then there is a sequence $a_1,\ldots, a_m \in A$ such that $y = f(x) = f(a_1)\cdots f(a_m)$, and so $B = \{f(a) : a \in A\}$ is a generating set for \mct.  Conversely, if $y \in \mct$ and $\ell_B(y) = m$, then there exist $a_1,\ldots,a_m\in A$ such that 
$y=f(a_1)\cdots f(a_m)$.  Let $x=a_1\cdots a_m \in \mcs$.  
Then $f(x) = y$ and $\ell_A(x) \leq m.$  This implies that $\gamma^{(T)}_B(n) \leq \gamma^{(S)}_A(n)$ for all $n$, and the growth conditions follow directly from this inequality.
\end{proof}

The growth of a finitely generated abelian semigroup  $\mathcal{S}$ is always polynomial.  Indeed, if $A$ is a set of generators for  $\mathcal{S}$ with $\card(A) = k$, then 
\[
\lambda_A(n) \leq {n+k-1 \choose k-1} \ll n^{k-1}
\]
and
\[
\gamma_A(n) \leq {n+k \choose k} \ll n^k.
\]
More precisely, Khovanskii~\cite{khov92,khov95}, Nathanson~\cite{nath00d}, Nathanson and Ruzsa~\cite{nath02a}  proved that there is a polynomial $f_A(x)$ with integer coefficients such that $\lambda_A(n) = f_A(n)$ for all sufficiently large integers $n$.  
It follows that there is a polynomial $F_A(x)$ with integer coefficients such that $\gamma_A(n) = F_A(n)$ for all sufficiently large $n$.  

The growth of finitely generated free semigroups of rank at least two is always exponential.  If \mcs\ is the free  semigroup generated by a set of $k \geq 2$ elements, 
then $\lambda_A(n) = k^n$ and $\gamma_A(n) = (k^{n+1}-1)/(k-1)> k^n.$
By Lemma~\ref{intermed:lemma:subsemi}, a  semigroup that contains a free subsemigroup on two generators has exponential growth.

Semigroups of intermediate growth are more difficult to construct.  
Beljaev, Sesekin, and Trofimov~\cite{belj-sese-trof77} proved that the free semigroup generated by two elements $e$ and $g$ with the relations $e^2 = e$ and $eg^ieg^je = eg^jeg^ie$ for all nonnegative integers $i$ and $j$ is a semigroup of intermediate growth.  
Okni{\` n}ski~\cite{okni98} constructed two examples of linear semigroups, that is, subsemigroups of the multiplicative semigroup $M_n(\Z)$ of $n \times n$ matrices.  Let $\mcs_1$ be the subsemigroup of $M_2(\Z)$ generated by the set
\[
A = \left\{ \left(\begin{matrix}
1 & 1 \\
0 & 1
\end{matrix}\right),
\left(\begin{matrix}
1 & 0 \\
1 & 0
\end{matrix}\right)
\right\}
\]
and let $\mct_1$ be the subsemigroup of $M_3(\Z)$ generated by the set
\[
A = \left\{ \left(\begin{matrix}
1 & 1 & 1\\
0 & 2 & 1\\
0 & 0 & 1
\end{matrix}\right),
\left(\begin{matrix}
1 & 0 & 0 \\
0 & 0 & 0 \\
0 & 0 & 1
\end{matrix}\right)
\right\}.
\]
Both $\mcs_1$ and $\mct_1$ are homomorphic images of the Beljaev, Sesekin, and Trofimov semigroup and have intermediate growth.   
Nathanson~\cite{nath99b} gave a simple number-theoretical proof that $\mathcal{S}^{(1)}$ has intermediate growth, and Lavrik-M\" annlin~\cite{lavr01} computed the growth function of $\mathcal{S}^{(1)}$.
Grigorchuk~\cite{grig88} and Shneerson~\cite{shne01a,shne04a,shne05a} have also investigated the growth of semigroups. 

In this paper we apply ideas from additive number theory to construct families of linear semigroups of intermediate growth that generalize Okni{\` n}ski's examples.

\section{Partition functions of  sets of positive integers}
Results about the intermediate growth of groups and semigroups often use estimates for the asymptotics of partition functions in additive number theory.  We review some of these results here.

Let $A$ be a set of positive integers and let $p_A(n)$ denote the number of partitions of $n$ into parts belonging to $A$.  If $A = \N$ is the set of all positive integers, then  $p(n) = p_{\N}(n)$ is the classical partition function.

If $A$ is a nonempty finite set of relatively prime positive integers and $\card(A) = r,$ then
\[
p_A(n) = \left( \frac{1}{\prod_{a\in A} a} \right)\frac{n^{r-1}}{(r-1)!} + O\left(n^{r-2}\right)
\]
(Nathanson~\cite{nath00a} and~\cite[Theorem 15.2]{nath00aa}).
Let $p_r(n)$ denote the number of partitions of $n$ into at most $r$ parts, let $\hat{p}_r(n)$ denote the number of partitions of $n$ into exactly $r$ parts, and let $\hat{q}_r(n)$ denote the number of partitions of $n$ into exactly $r$ distinct parts.  Since the number of partitions of $n$ into at most $r$ parts is equal to the number of partitions of $n$ into parts belonging to the set $A = \{1,2,\ldots,r\}$, it follows that
\[
p_r(n) = \frac{n^{r-1}}{r!(r-1)!} + O\left(n^{r-2}\right)
\]
and 
\beq  \label{intermed:pkn}
\hat{p}_r(n) = p_r(n)-p_{r-1}(n) = \frac{n^{r-1}}{r!(r-1)!} + O\left(n^{r-2}\right).
\eeq
We note that $p_r(n)$ and $\hat{p}_r(n)$ are increasing functions of $n$, and 
\beq      \label{intermed:prp}
p(n) =  \sum_{r=1}^n \hat{p}_r(n).
\eeq

If $n=a_1+\cdots + a_r$ is a partition of $n$ into $r$ distinct positive parts $a_1 > \cdots > a_r$, then $n \geq r(r+1)/2$ and 
\[
n - \frac{r(r-1)}{2} = (a_1- (r-1)) +( a_2 - (r-2)) + \cdots + (a_2 -1) + a_1
\]
is a partition of $n - r(r-1)/2$ into $r$ parts.  This identity establishes a bijection between partitions into exactly $r$ distinct parts and partitions into exactly $r$ parts, and so
\beq  \label{intermed:qkn}
\hat{q}_r(n) = \hat{p}_r(n - r(r-1)/2)  = \frac{n^{r-1}}{r!(r-1)!} + O\left(n^{r-2}\right).
\eeq

The set $A$ has asymptotic density  $d(A) = \alpha$ if
\[
\lim_{n\rightarrow\infty} \frac{1}{n} \sum_{\substack{a \in A\\ a\leq n}} 1 = \alpha.
\]
If $d(A) = \alpha>0$ and $\gcd(A) = 1$, 
then
\beq    \label{intermed:part}
\log p_A(n) \sim c_0\sqrt{\alpha n}
\eeq
where
\[
c_0 = 2\sqrt{\frac{\pi^2}{6}}.
\]
(An elementary proof of~\eqref{intermed:part} is in Nathanson~\cite[Theorem 16.1]{nath00aa}.)
In particular, if $A = \N$, then $\alpha = 1$ and the partition function $p(n) = p_{\N}(n)$ satisfies the Hardy-Ramanujan asymptotic estimate 
\[
\log p(n) \sim  c_0\sqrt{ n}.
\]
If $A$ is the set of odd positive integers, then $\alpha = 1/2$
and $\log p_A(n) \sim  c_0\sqrt{n/2}.$
Let $q(n)$ denote the number of partitions of $n$ into distinct positive integers.  Since the number of partitions of an integer into distinct parts is equal to the number of partitions into odd parts, it follows that
\[
\log q(n) \sim  c_0\sqrt{n/2}.
\]

\section{A condition for subexponential growth}

\bl       \label{intermed:lemma:gijcondition}
Let $\mathcal{S}$ be a semigroup generated by a set $A = \{e,g\}$, where $e$ is an idempotent.  The set 
\[
\mathcal{S}_0 = e\mathcal{S}e = \{exe:x\in \mathcal{S}\}
\]
is a subsemigroup of $\mathcal{S}$ with identity $e$, and is generated by the set $\{ eg^ke: k = 0,1,2,3,\ldots \}$.
Define $g^0=1.$  Then
\[
\mathcal{S} = \{g^i\}_{i=0}^{\infty} \cup \bigcup_{i,j=0}^{\infty} g^{i}\mathcal{S}_0g^{j}.
\]
Suppose that, for all nonnegative integers $i_1,j_1,i_2,j_2$,
\beq     \label{intermed:gijcondition}
g^{i_1}\mathcal{S}_0g^{j_1} \cap g^{i_2}\mathcal{S}_0g^{j_2} \neq \emptyset
\text{ if and only if } (i_1,j_1) = (i_2,j_2).   
\eeq
For all nonnegative integers $i$ and $j$, let 
 \[
 \gamma_{i,j}(n) = \card\{y\in g^i\mcs_0 g^j : \ell_A(y) \leq n\}
 \]
 and let $\gamma_0(n) = \gamma_{0,0}(n).$  Then
 \[
 \gamma_{i,j}(n) \leq \gamma_0(n-i-j).
 \]
 \el
 
\begin{proof}
Let $x \in g^i\mcs_0g^j.$  
By condition~\eqref{intermed:gijcondition}, every representation of $x$ as a word in $e$ and $g$ must be of the form $x = g^iy'g^j$ for some $y' \in \mcs_0$.  Therefore,
\[
\ell_A(x) = \min\{i+j+\ell_A(y') :   y' \in \mcs_0 \text{ and } x = g^iy'g^j \}
= i + j + \ell_A(y)
\]
for some $y \in \mcs_0$.
If $\ell_A(x) = m$, then $\ell_A(y) = m-i-j$  and so $\gamma_{i,j}(n) \leq \gamma_0(n-i-j).$
\end{proof}

\bt   \label{intermed:theorem:subexponential}
Let $\mathcal{S}$ be a semigroup generated by a set $A =
\{e,g\}$, where $e$ is an idempotent, and let $\mathcal{S}_0 =
e\mathcal{S}e.$ 
Suppose that
\benum
\item[(i)]   
For all nonnegative integers $k_1$ and $k_2$,
\[
eg^{k_1}eg^{k_2}e = eg^{k_2}eg^{k_1}e
\]
\item[(ii)]    
For all nonnegative integers $i_1,j_1,i_2,j_2$,
\[
g^{i_1}\mathcal{S}_0g^{j_1} \cap g^{i_2}\mathcal{S}_0g^{j_2} \neq \emptyset
\text{ if and only if } (i_1,j_1) = (i_2,j_2).
\]
\eenum
Then the semigroup $\mathcal{S}$ has subexponential growth.
\et

\begin{proof}
If $e=1$, then $\mcs = \mcs_0 = \{g^i\}_{i=0}^{\infty}$ and condition~(ii) is not satisfied.  Thus, we can assume that $e\neq 1$.
If $1 \in \mcs_0$, then  there exist positive integers $k_1,\ldots,k_r$ such that $eg^{k_1}eg^{k_2}e \cdots eg^{k_r}e=1.$  Multiplying this identity by $e$, we obtain $e=1$, which is absurd.  Therefore, $1 \notin \mcs_0$.

Let $\gamma_0(n)$ denote the number of elements  $y \in \mathcal{S}_0$ such that $\ell_A(y) \leq n.$
If $y \in \mathcal{S}_0$ and $\ell_A(y) = m \leq n$, then conditions~(i) and~(ii) imply that  either $y=e$ or there exist positive integers $r$ and $k_1 \geq \cdots \geq k_r$ such that
\[
y = eg^{k_1}eg^{k_2}e \cdots eg^{k_r}e
\]
and
\[
\ell_A(y) = k_1+k_2+\cdots + k_r+r+1 = m.
\]
Thus, to every $y \in \mcs_0$ with $\ell_A(y) =m$ there are associated  a positive integer $r$ and a partition of $m-r-1$ into exactly $r$ parts, and so the number of elements $y \in \mathcal{S}_0$ of length exactly $m$ is
\[
\lambda_0(m) = \gamma_0(m)-\gamma_0(m-1)
\leq \sum_{r=1}^m \hat{p}_r(m-r-1) \leq \sum_{r=1}^m \hat{p}_r(m) = p(m)
\]
by~\eqref{intermed:prp}.  Since $1 \notin \mcs_0$, we have
\[
\gamma_0(n)= \sum_{m=1}^n \lambda_0(m) \leq  \sum_{m=1}^n p(m) \leq np(n).
\]

Let $\gamma_{i,j}(n)$ denote the growth function of the set
$g^{i}\mathcal{S}_0g^{j}$ with respect to the generating set
$A$.   By Lemma~\ref{intermed:lemma:gijcondition},
\[
\gamma_{i,j}(n) \leq \gamma_0(n-i-j).
\]
By condition~(ii), if $g^i = g^j$, then $g^i\mcs_0 g^i = g^j \mcs_0 g^j$ and so $i = j.$
Therefore, for every $n \geq 0$ there is exactly one element in the set $\{g^i\}_{i=0}^{\infty}$ of length  $n$, and 
\begin{align*}
\gamma_A(n) & \leq \sum_{i=0}^{n} \sum_{j=0}^{n}\gamma_{i,j}(n) + n+1 \\
& \leq \sum_{i=0}^{n} \sum_{j=0}^{n}\gamma_{0}(n-i-j) + n +1\\
& \leq (n+1)^2\gamma_{0}(n) + n+1 \\
& \leq 2(n+1)^3 p(n).
\end{align*}
From the asymptotic formula~\eqref{intermed:part} for the partition function, we obtain
\[
\log \gamma_A(n) \leq \log 2(n+1)^3 + \log p(n) \ll \sqrt{n}
\]
and so
\[
\lim_{n\rightarrow\infty} \frac{\log \gamma_A(n)}{n} = 0.
\]
Thus, the growth function $\gamma_A(n)$ is subexponential.
\end{proof}

\section{Sequences with many partition products}
Let  $W=\{w_k\}_{k=1}^{\infty}$ be a sequence of elements in a semigroup $(\mathcal{X},\ast)$ such that
\[
w_{k_1}\ast w_{k_2} = w_{k_2} \ast w_{k_1}
\]
for all positive integers $k_1,k_2.$
To every finite sequence of positive integers $k_1, k_2, \ldots, k_r$ we associate the element $w_{k_1}\ast w_{k_2} \ast \cdots \ast w_{k_r} \in \mcx$.
For every integer $r \geq 1$, let $W_r(n)$ denote the subset of $\mathcal{X}$ associated to partitions of positive integers not exceeding $n$ into exactly $r$ parts, that is, $w \in W_r(n)$ if and only if there is a sequence of positive integers  $k_1 \geq k_2 \geq \cdots \geq k_r$ such that
$k_1 + k_2 + \cdots + k_r \leq n$ and $w = w_{k_1}\ast w_{k_2} \ast\cdots \ast w_{k_r}$.
We define
\[
\Phi(r) = \liminf_{n\rightarrow\infty} \frac{\log \card(W_r(n))}{\log n} .
\]
The sequence $\{w_k\}_{k=1}^{\infty}$ has \emph{many partition products in $W$} if
\[
\lim_{r\rightarrow\infty} \Phi(r) = \infty.
\]

Here are some examples of sequences with many partition products.  

\bt   \label{intermed:theorem:Ud}
For every positive integer $d$, let $\mcu_d$ be the multiplicative semigroup of all positive integers $u$ such that $u \equiv 1\pmod{d}.$
The sequence 
\[
W = \{dk+1\}_{k=1}^{\infty}
\]
has many partition products.
\et

\begin{proof}
Let $r \geq 1$ and let  $k_1 \geq k_2 \geq \cdots \geq k_r$ be a sequence of positive integers with $k_1+k_2+\cdots + k_r \leq n.$ 
Associated to this sequence is the integer $(dk_1+1)(dk_2+1)\cdots (dk_r+1) \in W_r(n).$
Let \mck\ denote the set of positive integers $k$ such that $dk+1$ is prime, and let $K(t)$ count the number of elements $k\in \mck$ with $k\leq t.$  
If $k_1 \geq k_2 \geq \cdots \geq k_r$ and $j_1 \geq j_2 \geq \cdots \geq j_r$ are distinct sequences of elements of \mck\ such that $k_1 \leq n/r$ and $j_1 \leq n/r$, then $k_1+k_2+ \cdots + k_r \leq n$ and $j_1+j_2+\cdots + j_r \leq n$, hence  $(dk_1+1)(dk_2+1)\cdots (dk_r+1) \in W_r(n)$ and $(dj_1+1)(dj_2+1)\cdots (dj_r+1) \in W_r(n)$.  By the fundamental theorem of arithmetic,  $(dk_1+1)(dk_2+1)\cdots (dk_r+1)  \neq (dj_1+1)(dj_2+1)\cdots (dj_r+1)$, and so
\[
\card(W_r(n)) \geq {K(n/r) + r-1 \choose r} \geq \frac{K(n/r)^r}{r!}.
\]
Let $\pi(x,d,1)$ denote the number of prime numbers $p \leq x$ such that $p\equiv 1\pmod{d}.$  Let
\[
x = \frac{dn}{r} + 1.
\]
Then $dk+1 \leq x$ is prime if and only if $k \leq n/r$.  
By the prime number theorem for arithmetic progressions, for sufficiently large $n$ we have 
\[
K(n/r) = \pi(x,d,1) \gg \frac{x}{\log x} \gg \frac{n}{r\log n}
\]
where the implied constants depend only on $d$.
Therefore,
\[
\card(W_r(n)) \geq \frac{K(n/r)^r}{r!} \gg \frac{n^r}{r!r^r(\log n)^r}
\]
and so
\[
\log\left( \card(W_r(n))\right) \gg  r\log n -\log (r!r^r) - r\log\log n
\]
and
\[
\Phi(r) = \liminf_{n\rightarrow\infty} \frac{ \log\left( \card(W_r(n))\right) }{\log n} \gg r.
\]
It follows that
\[
\lim_{r\rightarrow\infty} \Phi(r) = \infty
\]
and the sequence $\{dk+1\}_{k=1}^{\infty}$ has many partition products.  This completes the proof.
\end{proof}

In additive number theory, a sequence $\{b_k\}_{k=k_0}^{\infty}$ contained in an additive abelian semigroup is called a $\hat{B}_r$-sequence if sums of $r$ distinct terms of the sequence are distinct, that is, if $k_1 < k_2< \cdots < k_r$ and $j_1 < j_2 < \cdots < j_r$, and if
\[
b_{k_1} + b_{k_2} + \cdots + b_{k_r} = b_{j_1} + b_{j_2} + \cdots + b_{j_r}
\]
then $k_i = j_i$ for $i=1,\ldots,r.$  The sequence is called a $\hat{B}_{\infty}$-sequence if all finite sums of distinct elements of the set are distinct, that is, if $k_1<k_2< \cdots < k_r$ and $j_1 < j_2 < \cdots < j_s$, and if
\[
b_{k_1} + b_{k_2} + \cdots + b_{k_r} = b_{j_1} + b_{j_2} + \cdots + b_{j_s}
\]
then $r=s$ and $k_i = j_i$ for $i=1,\ldots,r.$    If $\{b_ k\}_{k=k_0}^{\infty}$ is a sequence of positive real numbers such that $\sum_{k=k_0}^{\ell-1} b_k < b_{\ell}$ for all $\ell > k_0$, then the sequence is a $\hat{B}_{\infty}$-sequence.  In particular, a $\hat{B}_{\infty}$-sequence is a $\hat{B}_r$-sequence for all $r \geq 1.$

\bl    \label{intermed:lemma:ctr}
Let $c,c_1,d,$ and $t$  be positive real numbers with $t > 2$, and let $\{e_k\}_{k=1}^{\infty}$ be a sequence of real numbers such that $|e_k| \leq c_1k^d$ for all $k \geq 1$.  There is an integer $k_0$ such that the sequence
\[
\{ ct^k + e_k\}_{k=k_0}^{\infty}
\]
is a strictly increasing $\hat{B}_{\infty}$-sequence.
\el

\begin{proof}
Let $b_k = ct^k + e_k$ for $k \geq 1.$
Since $t>2$ and $|e_k| \leq c_1k^d$,  there is an integer $k_0$ such that
\[
0 < b_k <  b_{k+1}
\]
for all $k \geq k_0$, and also
\beq      \label{intermed:partprodineq}
\frac{c_1(t-1)(k+1)^{d+1}}{c(d+1)(t-2)} < t^k
\eeq
for all $k > k_0.$  Let $\ell > k_0.$  
Using~\eqref{intermed:partprodineq} and the inequality
\[
\sum_{k=1}^{\ell} k^d < \frac{(\ell +1)^{d+1}}{d+1}
\]
we obtain
\begin{align*}
\sum_{k=k_0}^{{\ell}-1} b_k 
& = c\sum_{k=k_0}^{{\ell}-1}t^k +  \sum_{k=k_0}^{{\ell}-1}e_k  \\
& < \frac{ct^{\ell}}{t-1} + c_1 \sum_{k=k_0}^{{\ell}-1} k^d \\ 
& < ct^{\ell} - c_1{\ell}^d \\
& \leq b_{\ell}.
\end{align*}
This completes the proof.
\end{proof}

\bt      \label{intermed:theorem:Br}
Let $c,c_1d,$ and $t$  be positive real numbers with $t > 2$, and let $\{e_k\}_{k=1}^{\infty}$ be a sequence of real numbers such that $|e_k| \leq c_1k^d$ for all $k \geq 1$.  Let $k_0$ be a positive integer such that sequence
\[
W = \{ ct^k + e_k\}_{k=k_0}^{\infty}
\]
is a $\hat{B}_{\infty}$-sequence.  Then $W$ has many partial products.
\et

\begin{proof}
Let $b_k = ct^k + e_k$ for $k \geq 1$.   By Lemma~\ref{intermed:lemma:ctr}, there is an integer $k_0$ such that 
$W = \{ ct^k + e_k\}_{k=k_0}^{\infty}$ is a strictly increasing $\hat{B}_r$-sequence for every positive integer $r$.  
To every partition $n = k_1 +  \cdots + k_r$ into exactly $r$ distinct parts $ k_1 > k_2 > \cdots > k_r$, we associate the real number $b_{k_1+k_0} + b_{k_2+k_0} + \cdots + b_{k_r+k_0} \in  W_r(n+rk_0).$  Since $W$ is a $\hat{B}_r$-sequence, it follows that different partitions are associated to different real numbers, and so, by the partition asymptotic~\eqref{intermed:qkn},
\[
\card(W_r(n+rk_0)) \geq \hat{q}_r(n) \gg n^{r-1}
\]
and
\[
\Phi(r) = \liminf_{n\rightarrow\infty}\frac{\log \card(W_r(n))}{\log n} \geq r-1.
\]
This completes the proof.
\end{proof}

\section{A condition for superpolynomial growth}

\bt     \label{intermed:theorem:superpolynomial}
Let $\mathcal{S}$ be a semigroup generated by a set $A =
\{e,g\}$, where $e$ is an idempotent. Let $\mathcal{S}_0 =
e\mathcal{S}e = \{exe:x\in \mathcal{S}\}.$ Suppose that, for all
nonnegative integers $k_1$ and $k_2$,
\[
eg^{k_1}eg^{k_2}e = eg^{k_2}eg^{k_1}e.
\]
Let $\varphi$ be a semigroup homomorphism from $\mathcal{S}_{0}$
into a semigroup $(\mcx,\ast)$ such that the sequence
\[
W = \{\varphi(eg^k e)\}_{k=1}^{\infty}
\]
has many partition products in $\mcx$.  Then the semigroup $\mathcal{S}$ has superpolynomial growth.
\et

\begin{proof}
Let $w_k = \varphi(eg^ke)$ for $k = 1,2,\ldots.$
Then
\begin{align*}
w_{k_1}\ast w_{k_2} & = \varphi(eg^{k_1}e) \ast \varphi(eg^{k_2}e)
=  \varphi(eg^{k_1}eg^{k_2}e)  \\
& =  \varphi(eg^{k_2}eg^{k_1}e)
=  \varphi(eg^{k_2}e) \ast \varphi(eg^{k_1}e) \\
& = w_{k_2}\ast w_{k_1}
\end{align*}
for all $k_1, k_2 \in \N_0.$  If $y =  eg^{k_1}eg^{k_2}e \cdots eg^{k_r}e \in \mathcal{S}_0$, then
\[
\varphi(y) = \varphi\left(eg^{k_1}e\right)\ast
\varphi\left(eg^{k_2}e \right)\ast  \cdots
\ast\varphi\left(eg^{k_r}e\right) = w_{k_1}\ast  w_{k_2} \ast \cdots
\ast  w_{k_r}.
\]

Let $\gamma_A(n)$ denote the growth function of $\mathcal{S}$ with respect to $A$,
and let $\gamma_0(n)$ denote the growth function of $\mathcal{S}_0$ with respect to $A$.
Fix a positive integer $r$.  If $w \in W_r(n-r-1)$, then there are positive
integers $k_1, k_2, \ldots, k_r$ such that
\[
w_{k_1} \ast   w_{k_2} \ast   \cdots \ast   w_{k_r} = w
\]
and
\[
k_1 + \cdots + k_r \leq n-r-1.
\]
Let
\[
y = eg^{k_1}eg^{k_2}e\cdots eg^{k_h}e \in \mathcal{S}_0.
\]
Then
\[
\varphi(y) = w_{k_1} \ast   w_{k_2} \ast   \cdots \ast   w_{k_r} = w
\]
and
\[
\ell_A(y) \leq k_1+k_2+\cdots + k_r +r+1 \leq n.
\]
Thus, to every element $w \in W_r(n-r-1)$ there is at least one element $y\in \mathcal{S}_0$ with $\varphi(y) = w$ and $\ell_A(y) \leq n$, and so
\[
\gamma_0(n) \geq \card(W_r(n-r-1)).
\]
It follows that
\[
\liminf_{n\rightarrow\infty} \frac{\log\gamma_A(n)}{\log n}
\geq
\liminf_{n\rightarrow\infty} \frac{\log\gamma_0(n)}{\log n}
\geq \liminf_{n\rightarrow\infty} \frac{\log\card(W_r(n-r-1))}{\log n} =  \Phi(r).
\]
Since this inequality is true for all positive integers $r$ and $\Phi(r)$ tends to infinity, it follows that
\[
\lim_{n\rightarrow\infty} \frac{\gamma_A(n)}{\log n} = \infty
\]
and so the growth function $\gamma_A(n)$ is superpolynomial.
\end{proof}

\section{The semigroups $\mathcal{S}^{(d)}$}
\bt   \label{intemed:theorem:Sd}
Let $d$ be a positive integer and let $\mathcal{S}^{(d)}$ be the subsemigroup of $M_2(\Z)$ generated by the matrices
\[
g = \left(\begin{matrix}
1 & d \\
0 & 1
\end{matrix}\right)
\]
and
\[
e = \left(\begin{matrix}
1 & 0  \\
1 & 0
\end{matrix}\right).
\]
The semigroup $\mcs^{(d)}$ has intermediate growth.
\et

\begin{proof}
Multiplying matrices, we obtain
\[
eg^ke = (dk+1)e
\]
and 
\[
eg^{k_1}eg^{k_2}e = (dk_1+1)(dk_2+1) = eg^{k_2}eg^{k_1}e
\]
for all nonnegative integers $k_1,k_2.$  
It follows that $\mcs^{(d)}$ satisfies condition~(i) of 
Theorem~\ref{intermed:theorem:subexponential}, and 
\[
\mathcal{S}^{(d)}_0 = e\mathcal{S}^{(d)}e = \{ (dk+1)e : k = 0,1,2,\ldots\} = \{ue : u \in \mathcal{U}_d\}
\]
where $\mathcal{U}_d$ is the multiplicative semigroup of positive integers that are congruent to 1 modulo $d$.  

For all $u \in \mcu_d$ and all nonnegative integers $i$ and $j$,
\[
g^i ue g^j = 
u \left(\begin{matrix}
1 & di \\
0 & 1
\end{matrix}\right)
\left(\begin{matrix}
1 & 0  \\
1 & 0
\end{matrix}\right)
 \left(\begin{matrix}
1 & dj \\
0 & 1
\end{matrix}\right)
=
u\left(\begin{matrix}
di +1& dj+d^2 ij \\
1 & dj
\end{matrix}\right).
\]
Let $u_1,u_2 \in \mcu_d$ and $i_1,i_2,j_1,j_2 \in \N_0$.  Then
\[
g^{i_1} u_1e g^{j_1} = g^{i_2} u_2e g^{j_2}
\]
if and only if $u_1=u_2$ and $(i_1,j_1) = (i_2,j_2).$ 
Thus, the semigroup $\mcs^{(d)}$ also satisfies condition~(ii) of 
Theorem~\ref{intermed:theorem:subexponential}, and so has subexponential growth.

The function
\[
\varphi: \mathcal{S}_0^{(d)} \rightarrow \mathcal{U}_d
\]
defined by
\[
\varphi(ue) = u
\]
is a semigroup homomorphism.
Let $w_k = \varphi(eg^ke) = \varphi((dk+1)e) = dk+1$ and $W = \{dk+1\}_{k=1}^{\infty}$.  By Theorem~\ref{intermed:theorem:Ud}, the sequence $W$ has  many partition products.  By Theorem~\ref{intermed:theorem:superpolynomial}, the semigroup $\mathcal{S}^{(d)}$ has superpolynomial growth.
This completes the proof.
\end{proof}

\bp
If $d'$ divides $d$, then $\mathcal{S}^{(d)}$ is a subsemigroup of $\mathcal{S}^{(d')}$.  In particular, for every positive integer $d$, the semigroup $\mathcal{S}^{(d)}$ is a finitely generated subsemigroup of the Okni{\` n}ski semigroup $\mathcal{S}^{(1)}$.  
Let $\mathcal{S}$ be finitely generated subsemigroup of a finitely generated semigroup $\mathcal{S}$.  How are the growth rates of $\mathcal{S'}$ and $\mathcal{S}$ related?
\ep

\section{Subsemigroups of $M_3(\R)$}

For $k \geq 0$ and $m \geq 1$, let $F_k(x_1,\ldots,x_m)$ be the symmetric function of degree $k$ in $m$ variables defined by
\[
F_k(x_1,\ldots,x_m) = \sum_{\substack{(i_1,\ldots,i_m)\in \N_0^m \\ i_1+i_2+\cdots + i_m=k}} x_1^{i_1}\cdots x_m^{i_m}.
\]
Define  $F_k(x_1,\ldots,x_m)=0$
for $k<0$.   If $k \geq 1$ and $m \geq 2$, then 
\[
F_k(x_1,\ldots,x_m) = x_1F_{k-1}(x_1,x_2,\ldots,x_m) + F_k(x_2,\ldots,x_m).
\]
In particular,
\beq   \label{intermed:symmetric2}
F_k(x_1,x_2) = x_1F_{k-1}(x_1,x_2) + x_2^k
\eeq
and
\beq   \label{intermed:symmetric3}
F_k(x_1,x_2,x_3) = x_1F_{k-1}(x_1,x_2,x_3) + F_k(x_2,x_3).
\eeq

\bl     \label{intermed:lemma:semigroupT}
Consider the upper triangular matrices
\[
g = \left(\begin{matrix}
s & v & w \\
0 & t & u \\
0 & 0 & r
\end{matrix}\right)
\]
and
\[
e = \left(\begin{matrix}
1 & 0 & 0 \\
0 & 0 & 0 \\
0 & 0 & 1
\end{matrix}\right)
\]
with coefficients in a ring.  
For all $k \geq 0,$ define
\begin{align*}
u_k & = uF_{k-1}(r,t) \\
v_k & = vF_{k-1}(s,t) \\
w_k & = wF_{k-1}(r,s) + uv F_{k-2}(r,s,t).
\end{align*}
Then
\beq  \label{intermed:e2e}
e^2=e
\eeq
\beq  \label{intermed:matrixgk}
g^k = \left(\begin{matrix}
s^k & v_k & w_k \\
0 & t^k & u_k \\
0 & 0 & r^k
\end{matrix}\right)
\eeq
\beq  \label{intermed:matrix-egke}
eg^{k}e = \left(
\begin{matrix}
s^k & 0 & w_{k} \\
0 & 0& 0 \\
0 & 0 & r^k
\end{matrix}
\right)
\eeq
and
\beq  \label{intermed:matrix-egk1k2e}
eg^{k_1}eg^{k_2}e  = \left(\begin{matrix}
s^{k_1+k_2} & 0 & s^{k_1}w_{k_2}+r^{k_2}w_{k_1} \\
0 & 0& 0 \\
0 & 0 & r^{k_1+k_2}
\end{matrix}\right).
\eeq
For $k_1,k_2\in \N_0$, 
\beq  \label{intermed:k1k2commute}
eg^{k_1}eg^{k_2}e = eg^{k_2}eg^{k_1}e
\eeq
if and only if 
\beq  \label{intermed:k1k2commutebecause}
s^{k_1}w_{k_2}+r^{k_2}w_{k_1} = s^{k_2}w_{k_1}+r^{k_1}w_{k_2}.
\eeq
If $r=s=1$, then~\eqref{intermed:k1k2commute} holds for all $k_1,k_2\in \N_0.$
\el

\begin{proof}
Identity~\eqref{intermed:matrixgk} holds for $k=0$ and $k=1$.  If the formula is true for some 
$k \geq 1$, then~\eqref{intermed:symmetric2} and~\eqref{intermed:symmetric3} imply that 
\begin{align*}
g^{k+1} & = \left(\begin{matrix}
s^k & v_k & w_k \\
0 & t^k & u_k \\
0 & 0 & r^k
\end{matrix}\right)
 \left(\begin{matrix}
s & v & w \\
0 & t & u \\
0 & 0 & r
\end{matrix}\right) \\
& = \left(\begin{matrix}
s^{k+1} & s^kv+v_kt & s^kw+v_ku+w_kr\\
0 & t^{k+1} & t^ku+u_kr \\
0 & 0 & r^{k+1}
\end{matrix}\right) \\
& = \left(\begin{matrix}
s^{k+1} & v_{k+1} & w_{k+1}\\
0 & t^{k+1} & u_{k+1} \\
0 & 0 & r^{k+1}
\end{matrix}\right)
\end{align*}
and identity~\eqref{intermed:matrixgk} follows by induction on $k$.
Identities~\eqref{intermed:matrix-egke} and~\eqref{intermed:matrix-egk1k2e} follow by matrix multiplication, 
and~\eqref{intermed:k1k2commute} 
and~\eqref{intermed:k1k2commutebecause}
follow from inspection of~\eqref{intermed:matrix-egk1k2e}.
This completes the proof.
\end{proof}

\bt   \label{intermed:theorem:semigroupT}
Let $t,u,v,w$ be real numbers with $t > 2$.  Let $\mathcal{T}$ be the subsemigroup of $M_3(\R)$ generated by the matrices
\[
e = \left(\begin{matrix}
1 & 0 & 0 \\
0 & 0 & 0 \\
0 & 0 & 1
\end{matrix}\right)
\]
and
\[
g = \left(\begin{matrix}
1 & v & w \\
0 & t & u \\
0 & 0 & 1
\end{matrix}\right)
\]
The semigroup \mct\ has intermediate growth.
\et

\begin{proof}
For every $k\in \N$, let
\[
t_k = 1 + t + t^2 + \cdots + t^{k-1} = \frac{t^k-1}{t-1}
\]
and
\[
w_k = \frac{uv t^k}{(t-1)^2}  + \left(w-\frac{uv}{t-1}\right)k - \frac{uv}{(t-1)^2}.
\]
By Lemma~\ref{intermed:lemma:semigroupT},  we have $e^2=e$,
\[
g^k = \left(
\begin{matrix}
1 & vt_k & w_k \\
0 & t^k & ut_k \\
0 & 0 & 1
\end{matrix}
\right)
\]
\[
eg^ke = \left(
\begin{matrix}
1 & 0 & w_k \\
0 & 0 & 0 \\
0 & 0 & 1
\end{matrix}
\right)
\]
and
\[
eg^{k_1}eg^{k_2}e = eg^{k_2}eg^{k_1}e
\]
for $k_1,k_2\in \N_0$.

Let $\mct_0 = e\mct e.$  Then $\mct_0$ consists of all matrices of the form 
\[
eg^{k_1}eg^{k_2}e\cdots eg^{k_r}e = 
\left(
\begin{matrix}
1 & 0 & w_{k_1}+w_{k_2}+\cdots + w_{k_r} \\
0 & 0 & 0 \\
0 & 0 & 1
\end{matrix}
\right)
\]
where $k_1,\ldots, k_r$ is a finite sequence of nonnegative integers.
For $i,j \in \N_0$, the set $g^i\mct_0 g^j$ consists of all matrices of the form
\beq   \label{intermed:eMatrix}
g^ieg^{k_1}eg^{k_2}e\cdots eg^{k_r}e g^j= 
\left(
\begin{matrix}
1 &  vt_j & w_i + w_j + w_{k_1}+w_{k_2}+\cdots + w_{k_r} \\
0 & 0 & ut_i \\
0 & 0 & 1
\end{matrix}
\right)
\eeq
Since $\{t_i\}_{i=0}^{\infty}$ is a strictly increasing sequence of real numbers, 
formula~\eqref{intermed:eMatrix} implies that for all $i_1,i_2,j_1,j_2 \in \N_0$, 
\[
g^{i_1}\mct_0 g^{j_1} \cap g^{i_2}\mct_0 g^{j_2} \neq \emptyset
\]
if and only if 
\[
(i_1,j_1) = (i_2,j_2).
\]
Thus, the semigroup \mct\ satisfies conditions~(i) and~(ii) of Theorem~\ref{intermed:theorem:subexponential} and so has subexponential growth.  

Define the function $\varphi:\mct_0 \rightarrow \R$ by
\[
\varphi(eg^ke) = w_k = ct^k + e_k
\]
where $c = uv/(t-1)^2$ and $e_k = (w-uv/(t-1))k -c = O(k)$. Then $\varphi$ is a homomorphism from $\mct_0$ into the additive group of real numbers.  By Theorem~\ref{intermed:theorem:Br}, there is an integer $k_0$ such that the sequence $W = \{w_k\}_{k=k_0}^{\infty}$ has many partial products.  By Theorem~\ref{intermed:theorem:superpolynomial}, the semigroup \mct\ has superpolynomial growth.
This completes the proof.
\end{proof}

\emph{Acknowledgements.}
I wish to thank David Newman, Kevin O'Bryant, and Lev Shneerson for helpful discussions about this work.

\def\cprime{$'$} \def\cprime{$'$} \def\cprime{$'$}
\providecommand{\bysame}{\leavevmode\hbox to3em{\hrulefill}\thinspace}
\providecommand{\MR}{\relax\ifhmode\unskip\space\fi MR }
\providecommand{\MRhref}[2]{%
  \href{http://www.ams.org/mathscinet-getitem?mr=#1}{#2}
}
\providecommand{\href}[2]{#2}

\end{document}